\documentclass[reqno,12pt]{amsart}

\usepackage{geometry,txfonts}
\usepackage[all]{xy}

\numberwithin{equation}{section}

\theoremstyle{plain}
\newtheorem{lemma}{Lemma}[section]
\newtheorem{teo}[lemma]{Theorem}

\newtheorem{propo}[lemma]{Proposition}

\theoremstyle{definition}
\newtheorem{defi}[lemma]{Definition}
\newtheorem{defs}[lemma]{Definitions}

\newtheorem{obs}[lemma]{Observation}

\theoremstyle{remark}

\newcommand{\pic}{{\rm Pic\thinspace}}
\newcommand{\bs}{{\rm Bs\thinspace}}

\newcommand{\p}{\mathbb{P}}
\newcommand{\pp}{\mathbb{P}}

\newcommand{\zz}{\mathbb{Z}}
\newcommand{\nat}{\mathbb{N}}

\newcommand{\oc}{{\mathcal O}}
\newcommand{\Osh}{{\mathcal O}}

\newcommand{\nn}{{\mathcal N}}

\newcommand{\aut}{\operatorname{Aut}}

\newcommand{\map}{\rightarrow}

\newcommand{\rt}{\longrightarrow}
\newcommand{\rmap}{\dashrightarrow}

\begin{document}

\setlength{\baselineskip}{15pt}

\title[On multiples of divisors associated to Veronese embeddings]{On multiples of divisors associated to Veronese embeddings with defective secant variety}
\author{Antonio Laface}
\address{
Departamento de Matem\'atica \newline
Universidad de Concepci\'on \newline 
Casilla 160-C \newline
Concepci\'on, Chile}
\email{alaface@udec.cl}
\email{}

\author{Luca Ugaglia}
\address{
Via Natta 2 \newline
14039, Tonco, Italia}
\email{luca.ugaglia@gmail.com}

\begin{abstract}
In this note we consider multiples $aD$, where
$D$ is a divisor of the blow-up of $\pp^n$ along points in general position which appears in the Alexander and Hirschowitz list of Veronese embeddings having defective secant varieties.
In particular we show that there is such a $D$ with $h^1(X,D) > 0$ and  $h^1(X,2D) = 0$.
\end{abstract}
\maketitle

\section*{Introduction}
In this note we deal with the problem of determining the cohomology groups of line bundles on the blow-up of the projective space $\pp^n$ at points in general position.
This problem is strictly related to that of finding the dimension of the $(r-1)$-secant variety, $S_{r-1}Z$, of the image of $d$-Veronese embedding: $\nu_{n,d}:\pp^n\map Z\subset\pp^N$. The two problems are connected by Terracini's Lemma~\cite{tev}.
In fact this lemma allows one to determine the dimension of $S_{r-1}Z$ once one knows $h^0(X_r^n,D)$, where  $D=\pi_r^*~\Osh_{\p^n}(d)- 2E_1-\dots-2E_r$ and $\pi_r: X_r^n\map \pp^n$ is the blow-up of the projective space along $r$ points in general
position with exceptional divisors $E_1,\dots,E_r$. 

It is standard to say that $D$ is {\em special} if $h^0(X_r^n,D)\cdot h^1(X_r^n,D) > 0$. In this case $h^0(X_r^n,D)$ is strictly bigger than expected. The problem of finding all the special $D$ of this form has been completely solved by  Alexander and Hirschowitz (see~\cite{ah} and~\cite{cha}) and the special divisors are listed in Table 1. 
A natural generalization of this problem is to study the special divisors of the form $D=\pi_r^*~\Osh_{\p^n}(d)-m_1E_1-\dots -m_rE_r$, with $m_i\in\nat$. In the case $n=2$ the solution to the problem  is predicted by the equivalent conjectures of Segre, Harbourne, Gimigliano and Hirschowitz
\cite{seg,har,gim,hi}, hereafter ``SHGH conjecture", that can be formulated in the following way: the divisor $D$ is special if and only if there exists a rational curve $C$ whose normal bundle is $\Osh_{\p^1}(-1)$ and such that $D\cdot C\leq -2$. 

When $n\geq 3$ it is true that, under some extra-hypothesis,  if there is a rational curve $C$ whose normal bundle is $\Osh_{\p^1}(-1)^{n-1}$ and such that $D\cdot C \leq -2$, then $D$ is special (see Theorem~\ref{h1}). The converse of this statement turns out to be false, since for instance in $\p^3$ there exist special divisors having non-negative intersection product with any rational curve. 
Anyway, in all the known cases of special divisors in $\p^2$ and $\p^3$, the {\em stable base locus} of $D$ (i.e. the intersection of the base loci of all the multiples of $D$) turns out to be at least $1$-dimensional and, as a consecuence of this, every multiple $aD$ of a special divisor is still special.
When $n\geq 4$ it seems that the situation can be much more complicated, since for instance there are examples of special divisors $D$ having stable base locus which is at most $0$-dimensional and such that there exist multiples of $D$ that are non-special.

The aim of this note is to study the behavior of the multiples of divisors appearing in Table 1, since they turn out to give examples of all the different phenomena we presented above. 
The note is structured as follows.  In the first section we fix the notations and we recall the list of Alexander 
and Hirschowitz and in Section 2 we prove Theorem~\ref{h1}. 
Then we begin an analysis of all the examples of the special divisors 
of Table 1 and their multiples; in Section 3 and 4 we consider the cases that can be explained by means of Theorem~\ref{h1}, while in Section 5 we study an example of a special divisor having non negative intersection with any rational curve and Section 6 deals with the last example, i.e. a system 
whose stable base locus  is at most $0$-dimensional and having a multiple which is not special. Finally, in the last section we propose some open problems related to the topics of the note.

\section{Preliminaries}

In what follows we will adopt the following notation:

\begin{itemize}\itemsep 5mm \leftskip -5mm
\item[-] $\pi_r: X_r^n\map\pp^n$ will be the blow-up of $r$ points in general position. The exceptional divisors will be denoted by $E_1,\ldots,E_r$ and the pull-back $\pi_r^*\Osh_{\pp^n}(1)$ by $H$.
\item[-] We denote by $D\boxtimes D'$ the divisor defined on $X\times X'$ as $p_1^*(D)\otimes p_2^*(D')$, where the $p_i$'s are the  projections on the two factors.
\item[-] Given a divisor $D$ of $X$, its {\em base locus} $\bs |D|$ is the intersection 
of all the $D'\in|D|$, while its {\em stable base locus} is the intersection of $\bs |nD|$ for all positive $n\in\zz$. 

\end{itemize}

\begin{defs}
Given a divisor $D:=dH-m_1E_1-\dots -m_rE_r$ we put
\[
v(D)=\binom{d+n}{n}-\sum_{i=1}^r\binom{n+m_i-1}{n}-1
\]
and we call $v(D)$ the { \em virtual dimension of } $D$ .
We denote by $e(D)$ the {\em expected dimension of} $D$, defined as follows:
\[
e(D)=max\{-1,v(D)\}.
\]
We say that $D$ is {\em special} if $\dim|D|>e(D)$ and
its {\em speciality} is given by 
\[
\dim|D|-v(D).
\]
\end{defs}

\begin{defi}
Given a birational map $\phi: X\rmap Y$, we say that $\phi$ is {\em small} if there are Zariski open subsets $U, V$, such that $\phi_{U}: U\map V$ is an isomorphism and the codimension of $X\setminus U$ and $Y\setminus V$ is at least $2$. Any small birational map induces a pull-back isomorphism $\phi^*: \pic(Y)\map \pic(X)$ defined by $\phi^*D := \overline{\phi^{-1}(D\cap V)}$.
\end{defi}

Let us consider now the following lemma which allows us to compare the cohomology of line bundles and their pull-back in a blow-up.

\begin{lemma}\label{compare}
Let $\pi: Y\map X$ be a blow-up map of a smooth algebraic variety $X$ along a subvariety $V$. For any divisor $D$ of $X$ and any non-negative $i$, we have that:
\[
h^i(X,D) = h^i(Y,\pi^*D).
\]
\end{lemma}
\begin{proof}
By Zariski connectedness theorem we have that $\pi_*\Osh_Y=\Osh_X$ (see~\cite[Chapter III, \S 11]{hart}). This together with the projection formula (see~\cite[Chapter II, \S 5]{hart}) implies:
\begin{equation}\label{proj}
\pi_*\pi^*D = \pi_*(\Osh_Y\otimes \pi^*D) \cong
\pi_*\oc_Y\otimes D = D.
\end{equation}
Again by the projection formula (see~\cite[Chapter III, \S 8]{hart}) we also have that $R^i\pi_*\pi^*D = 0$ for any $i\geq 1$. Therefore $H^i(Y,\pi^*D)\cong
H^i(X,\pi_*\pi^*D)$, and by equation~\ref{proj} we get the
thesis.

\end{proof}

In this note we will focus on linear systems of the form $|aD|$ where
\[
D = dH-2E_1-\dots-2E_r
\]
and $(d,n,r)$ belong to the following table (see for instance~\cite{bo}):

\begin{center}
\begin{table}[h]
\caption{}
\begin{tabular}{c|c|c|c}
degree $d$ & dimension $n$ & n. of points $r$ & degree of $C$ \\
\hline
$2$ & $n\geq 2$ & $2\leq r\leq n$ & 1 \\
$3$ & $4$ & $7$ & $3$ \\
$4$ & $2$ & $5$ & $2$ \\
$4$ & $3$ & $9$ & \\
$4$ & $4$ & $14$ & 
\end{tabular}
\end{table}
\end{center}

The reason for the fourth column of the table will be explained in the next section.

\section{Speciality and rational curves}

Let $C\subset X$ be a smooth rational curve whose normal bundle $\nn_{C|X}$ is a direct sum of copies of $\Osh_{\p^1}(-1)$.

\begin{teo}\label{h1}
Let $D$ be an effective divisor such that $D\cdot C = - t < 0$.
Let $\pi: Y\map X$ be the blow-up along $C$ with exceptional divisor $E$, then:
\[
h^1(X,D) = h^1(Y,\pi^*D-tE) - h^2(Y,\pi^*D-tE) +
\binom{t+\dim X-2}{t-2}.
\]
\end{teo}
\begin{proof}
Consider the exact sequence of sheaves:
\[
\xymatrix{
0\ar[r] & \pi^*D-(r+1)E\ar[r] & \pi^*D-rE\ar[r] & \pi^*D-rE_{|E}\ar[r] & 0.
}
\]
Summing-up the Euler characteristics of the preceding sequence as $r$ varies from $0$ up to $t-1$, by Lemma~\ref{compare} we obtain:
\[
\chi(X,D) = \chi(Y,\pi^*D-tE) + \sum_{r=0}^{t-1}\chi(E,\pi^*D-rE_{|E}).
\]
The exceptional divisor $E$ is isomorphic to $\pp^1\times\pp^{n-2}$, where $n$ is the dimension of $X$. Denoting by $H_i$ the hyperplane class of $\pp^i$ we have $E_{|E} = -H_1\boxtimes -H_{n-2}$ so that: 
\[
\pi^*D-rE_{|E} = (r-t)H_1\boxtimes rH_{n-2}.
\]
Observe that this sheaf has no non-zero global sections if $r < t$. Looking at the exact sequence this means that
$h^0(X,D) = h^0(Y,\pi^*D-tE)$.
On the other hand, the Kunneth formula (see for instance~\cite{kun}) gives us
\[
h^1(E,(r-t)H_1\boxtimes rH_{n-2}) = (t-r-1)\binom{n-2+r}{r}
\]
and $h^i(E,(r-t)H_1\boxtimes rH_{n-2}) = 0$ for $i\geq 2$. This in turn gives the equalities
\[
h^i(Y,\pi^*D) = h^i(Y,\pi^*D-tE)
\]
for $i\geq 3$ and the thesis follows.
\end{proof}

Theorem~\ref{h1} suggests the idea that $D$ can have a positive $h^1$ in case $D\cdot C \leq -2$ for a rational $C$ as in the hypothesis. This turns out to be true in all the known examples.

\section{Lines and conics in the base locus of $|aD|$}

The $(2,n,r)$ and $(4,2,5)$ cases are well known in the literature. 
We start by considering the divisor $D:=H-E_1-\dots -E_r$, corresponding to the hyperplanes of $\p^n$ through $r$ general points $p_1,\dots, p_r$. We are going to study the systems $|dD|$, for $d\geq 2$ (these include the $(2,n,r)$ case and all its multiples).
For any $2\leq r\leq n$, the elements of the system $|dD|$ are cones whose vertex $V$ is the $\p^{r-1}$ spanned by the $r$ points. The projection
$\pi_V:\pp^n\rmap\p^{n-r}$ gives a bijection between the elements of $|dD|$ and those of $|\Osh_{\pp^{n-r}}(d)|$. Therefore the system is special and its dimension is
\[
\binom{n-r+d}{d}-1.
\]
Let us spend a few more words about these systems in some particular cases. We denote by $C_{ij}$ the strict transform of the line through the points $p_i$ and $p_j$ and by $E_{ij}$ the exceptional divisor corresponding to the blow up $\pi:Y\rt X$ along the $C_{ij}$'s.\\
First of all we remark that when $r=2$, the speciality of the system is
\[
\binom{d+n-2}{d-2},
\]
as predicted by Theorem~\ref{h1}, since $D\cdot C_{12}=-d$. Therefore in this case $h^2(Y,\pi^*(dD)-dE_{12})=0$.
Moreover, if $r$ is bigger than $2$, but $d=2$, the speciality is 
\[
\frac{r(r-1)}{2}.
\]
By Theorem~\ref{h1}, each line $C_{ij}$ gives a contribution 
of $1$, and hence, since there are $r(r-1)/2$ of them, $h^2(Y,\pi^*(2D)-\sum 2E_{ij})=0$ as before.

Let us see that when both $d$ and $r$ are bigger than $2$, the $h^2$ can be positive. 
For instance, if we consider 
the system corresponding to $3D:=3(H-E_1-E_2-E_3)$, 
we have that its virtual and effective dimensions are $-11$ and $0$ respectively.
But by Theorem~\ref{h1}, each of the $3$ lines gives a contribution of $4$, and this implies that 
in this case $h^2(Y,\pi^*(3D)-\sum 3E_{ij})=1$ (see~\cite[Example 7.6]{lu} for a more detailed description).

\vspace{.5cm}

The $(4,2,5)$ case is also well known. The only element of the linear system is twice the $(-1)$-curve $C := 2H-E_1-\dots-E_5$. 
Every multiple of the same $(-1)$-curve gives the special system $|aC|$, whose speciality is $a-1$. Observe that there is an automorphism $\phi\in\aut(X_5^2)$, induced by a quadratic transformation of $\pp^2$ such that $\phi(C)$ is of the form $H-E_1-E_2$, so that $2\phi(C)$ is of type $(2,2,2)$ which we already considered before.

\vspace{.5cm}

\section{The rational normal curve of $\pp^4$ and a small birational map.}

This case is similar to those of the preceding section because $|D|$ is special because of a rational curve $C$ with $C\cdot D \leq -2$. In this case the dimension of $|aD|$ can be determined by transforming the linear system by means of a small birational map.

Let $p_1,\ldots,p_5\in\pp^4$ be the points where all but one of their coordinates vanish and let $p_6$ and $p_7$ be two points in general position with respect to the first five.
Consider the Cremona transformation of $\pp^4$ given by 
\[
\sigma([x_0:\ldots:x_4]) := [x_0^{-1}:\ldots:x_4^{-1}].
\]
The first five points belong to the indeterminacy locus of $\sigma$.
The map $\sigma$ induces a small birational map $\phi: X_7^4\rmap {X'}_7^4$ between
two blow-ups of $\pp^4$ along seven points in general position.
The isomorphism $\phi^*: \pic({X'}_7^4)\map\pic(X_7^4)$ is described by the following.
\begin{propo}\label{cre-surfaces}
Let $D = dH-m_1E_1-\dots-m_7E_7$, then
\[
\phi^*D  =  D+k(H-E_1-\ldots-E_5),
\]
where $k=3d-m_1-\dots-m_5$.
\end{propo}

\begin{proof}
This is done by evaluating the action of $\sigma$ on monomials.
Let $T\subset\zz^5$ be the set of exponents of monomials of the linear system of hypersurfaces of degree $d$ with points of multiplicities $m_1,\ldots,m_5$ at $p_1,\ldots,p_5$ which is defined by:
\[
T := \{(a_1,\ldots,a_5)\in\zz^5\ |\ 0\leq a_i\leq d-m_i\text{ and } a_1+\dots+a_5=d\}.
\]
The action of $\sigma$ on the exponents of monomials is given by the involution:
\[
\sigma_{\zz}(a_1,\ldots, a_5) := (d-a_1-m_1,\ldots,d-a_5-m_5).
\]
Observe that if $(a_1,\ldots,a_5)\in T$ then $0\leq d-a_i-m_i\leq (d+k) - (m_i+k)$ and $\sum_{i=1}^5 (d-a_i-m_i) = d + k$, so that $\sigma_{\zz}(T)$ corresponds to the set $T_k$ of exponents of monomials of hypersurfaces of degree $d+k$ with multiplicities 
$m_1+k,\ldots,m_5+k$. Since $\sigma_{\zz}$ is an involution, this implies that $\sigma_{\zz}(T_k)\subseteq T$ as well and we get the thesis.
\end{proof}

Consider now the divisor 
\[
D = 3H-2E_1-\dots-2E_7.
\]
The linear system $|D|$ has virtual dimension $-1$ and it
contains only one element, namely the strict transform of the secant variety of the rational normal curve of $\p^4$, which is singular exactly along this curve (see~\cite{ch}).

Another way to prove the same result is to apply two transformations like $\phi$ and based on the first $5$ points and the last $5$ respectively, to $D$. We get
\begin{equation}\label{cubic}
\begin{array}{ccl}
D & \map & 2H-E_1-\dots-E_5-2E_6-2E_7 \\[5pt]
& \map & H-E_1-\dots-E_4,
\end{array}
\end{equation}
and in particular $\dim|D|=\dim |H-E_1-\dots-E_4|=0$.
This idea can be used to determine the dimension of $|aD|$ for any $a\geq 1$. The virtual dimension of this system is $-\binom{a+2}{2}$, while it is not empty, because it contains at least $aD$. Let us see
that in fact it does not contain any other divisor.
If we apply the same elementary transformations as before, we obtain
that 
\[
\dim |aD| = \dim |a(H-E_1-\dots-E_4)| = 0,
\]
since this system consists of $a$ times the hyperplane through the four multiple points. 

Observe that the strict transform $C$ of the rational normal curve through $p_1,\ldots,p_7$ has intersection $C\cdot D = -2$.
In this case we have that the speciality of $|aD|$ is exactly the one predicted by Theorem~\ref{h1}.

\section{Multiples of the strict transform of a quadric of $\pp^3$}

Here we consider the divisor given by the strict transform of a smooth quadric $Q$ through $9$ points in general position:
\[
D = 2H-E_1-\dots-E_9.
\]
The expected dimension of $|aD|$ is $-1$ and we are going to see that in fact $|aD| = aD$. By the adjunction formula, the restriction $D_{|D}$ coincides with the anticanonical bundle $-K_D$. 
Observe that an irreducible element $B\in |-aK_D|$ would have genus
$(a-a^2)/2+1$, so that $a\leq 2$. But this can not happen since otherwise the $9$ points would be in special position on the quadric.
Therefore the divisor $D$ is contained in the base locus of $|aD|$ for any $a\geq 1$.

We are now going to see a great difference between this case and the preceding ones. 
In fact, we can explain all the preceding cases by means of Theorem~\ref{h1}, since in each case there is a rational curve $C\subset X_r^n$ having negative intersection with $D$.
The following theorem shows that this is no longer the case.

\begin{teo}\label{curve}
Let $C$ be an irreducible rational curve of $X_9^3$ and let $D := 2H-E_1-\dots-E_9$, then $C\cdot D\geq 0$.
\end{teo}
\begin{proof}
Suppose that there exists such a $C$ with $C\cdot D < 0$. 
Consider the blow-up $Y_9$ of $\pp^3$ along nine points lying on a smooth elliptic quartic curve $\Gamma\subset\pp^3$. There is a flat family whose general fiber is isomorphic to $X_9^3$ and whose central fiber is isomorphic to $Y_9$. In this degeneration, the linear system  $|D|$ is degenerated into the strict transform $|D'|$ of a pencil of quadrics, while $C$ is sent to a, possibly reducible, new rational curve $C'$. Since $C'\cdot D' = C\cdot D < 0$, there is an irreducible component $C''$ of $C'$ which has negative intersection with $D'$.
Observe that $C''$ is a rational curve which must be contained in the base locus of $|D'|$. Since this base locus is an irreducible elliptic curve we get a contradiction.
\end{proof}

In fact, we can prove that this system is special (see~\cite[Conjecture 4.1]{lu}) by
observing that $h^1(X,2D) \geq h^1(D,2D_{|D})$. By the Riemann-Roch formula we see that $\chi(D,2D_{|D}) = -2$ so that the last cohomology group must have positive dimension.

\vspace{.5cm}

\section{Multiples of the strict transform of a quadric of $\pp^4$}

The aim of this section is to show, with a bit of computer algebra help, that there is a special linear system $|2D|$ on the blow-up of $\pp^n$ along points in general position, whose stable base locus is at most $0$-dimensional.
As in the preceding section we consider the divisor
\[
D := 2H-E_1-\dots-E_{14}
\]
defined on the blow-up of $\pp^4$ along $14$ points in general position.

Let us show a first different phenomenon respect to the preceding cases.

\begin{propo}
The divisor $D$ is nef and big.
\end{propo}
\begin{proof}
We proceed as for the proof of Theorem~\ref{curve} by specializing the $14$ points on the intersection $Q_1\cap\dots\cap Q_4$ of four general quadrics. Denote by $D'$ the divisor obtained by blowing-up the points in the new configuration and observe that the 
base locus of $|D'|$ is zero dimensional.
Assume that there exists an integral curve $C$ with $C\cdot D < 0$. After specializing the points, we have that $C$ deforms into a curve $C'$ which has a component $C''$ such that $C''\cdot D' < 0$. This in particular implies that $C''$ must be contained in the base locus of $|D'|$ which is a contradiction.

The bigness of $D$ comes from the fact that $D^4=2$.

\end{proof}

We conclude by proving that, unlike the previous cases, there is a multiple of $D$ which enjoys the following properties: 
\begin{propo}
$4D$ is non-special and the dimension of $\bs |4D|$ is at most $0$.
\end{propo}
\begin{proof}
Assume that the coordinates of the $p_i$'s are rational numbers, then there is a primitive integer vector representing each point in the projective space. Let $I$ be the homogeneous ideal of the $14$ quadruple points of $\pp^4$. We specialize this ideal as $J = I \pmod{101}$ and we evaluate $J$ by means of the following {\em Singular}~\cite{GPS05} program:
{\small
\begin{verbatim}
ring r = 101,x(0..4),dp;
int i; ideal J = 1; intmat a = random(3,14,4);
for (i=1;i<=14;i=i+1)
{
  ideal J(i) = 
  x(1)-a[i,1]*x(0),
  x(2)-a[i,2]*x(0),
  x(3)-a[i,3]*x(0),
  x(4)-a[i,4]*x(0);
 J = std(intersect(J(i)^4));
}
\end{verbatim}
}
This gives $\dim J_8 \leq 5$ for the degree $8$ homogeneous part of $J$. We deduce that $|4D|$ is non special since its expected dimension is $4$.
The following
{\small
\begin{verbatim}
Ideal B = J[1..5]; B = std(B);
\end{verbatim}
}
gives the ideal $B$ of the base locus of the linear system generated by the elements of $J_8$ which turns out to be zero dimensional.
\end{proof}

\section{Conclusions}

We have studied the cohomology of divisors $D$ on the blow-up $\pi_r: X_r^n\map\pp^n$ along points in general position. Our aim was to explore examples of effective divisors $D$ with $h^1(X_r^n,D)> 0$ which we called special.
To achieve this we moved our steps from the well known case of divisors of the form
\[
D = \pi_r^*\Osh_{\pp^n}(d)-2E_1-\dots -2E_r
\]
and concentrated our attention on its multiples $aD$. 
A natural observation which can be deduced from our work is:
\begin{obs}\label{bs}
In all the examples of this note, if $D$ is special then $\dim\bs |D| > 0$. Moreover, if $n\leq 3$ then $\dim \bigcap_{a\in\zz} \bs |aD| > 0$.
\end{obs}
We found that $\bs|D|$ was either a rational curve $C$ (Sections 2,3,4) or the blow-up of a smooth quadric (Sections 5 and 6). About the former case, in Theorem~\ref{h1} we proved a general statement on the speciality of divisors due to rational curves. 

Our work thus gives a picture of the main open problem left in this area:

{\bf Problem:} Prove that for any special $D\subset X_r^n$ there is an irreducible and reduced positive-dimensional variety $B\subseteq\bs|D|$ such that $\chi(B,D_{|B}) < 0$.

To the authors knowledge, in all the known cases $B$ is a rational variety. This is no longer true if one works with blow-ups of other varieties at points (see~\cite{dl} for an example with $K3$ surfaces). The SHGH conjecture predicts that, if $n=2$, then $B$ is a $(-1)$-curve. In~\cite{lu} is it conjectured that if $n=3$, then $B$ is either a rational curve or a surface which is transformed by a small birational map of $X_r^3$ into the blow-up of a smooth quadric. 
Recently, in~\cite{du}, the author introduced a new technique for determining the speciality of divisors of $X_r^n$. Perhaps this can be used for formulating a conjecture for any dimension $n$.

We leave the reader with the following observation which makes use of the small birational maps of $X_r^n$.

\begin{obs}
In all the examples that we considered with $D$ special, there is a small birational map $\phi$ of $X_r^n$ such that $\phi(B)$ is the blow-up of a smooth Fano variety.
\end{obs}

\end{document}